%% file: main.tex
\documentclass[english,conference]{ieeeconf}
\usepackage[T1]{fontenc}
\usepackage[latin9]{inputenc}
\usepackage{verbatim}
\usepackage{mathtools}
\usepackage{amsmath}
\usepackage{amsthm}
\usepackage{amssymb}
\usepackage{graphicx}

\makeatletter
\theoremstyle{plain}
\newtheorem{thm}{\protect\theoremname}
\theoremstyle{plain}
\newtheorem{prop}[thm]{\protect\propositionname}
\theoremstyle{remark}
\newtheorem{rem}[thm]{\protect\remarkname}
\theoremstyle{definition}
\newtheorem{defn}[thm]{\protect\definitionname}
\theoremstyle{plain}
\newtheorem{lem}[thm]{\protect\lemmaname}
\theoremstyle{plain}
\newtheorem{cor}[thm]{\protect\corollaryname}

\input{preamble.tex}
\input{affiliation.tex}
\makeatother

\usepackage{babel}
\providecommand{\corollaryname}{Corollary}
\providecommand{\definitionname}{Definition}
\providecommand{\lemmaname}{Lemma}
\providecommand{\propositionname}{Proposition}
\providecommand{\remarkname}{Remark}
\providecommand{\theoremname}{Theorem}

\begin{document}

\title{\textbf{\Large{}Gradient Methods with Dynamic Inexact Oracles}}

\author{Shuo Han\allthanks}
\maketitle
\begin{abstract}
\input{abstract.tex}
\end{abstract}

\section{Introduction}

\input{intro.tex}

\section{Mathematical Preliminaries}

\input{prelim.tex}

\section{Dynamic Inexact Oracles}

\input{oracle.tex}

\section{Convergence Analysis\label{sec:analysis}}

\input{analysis.tex}

\section{Conclusions}

\input{conclusions.tex}

\appendix{}

\input{appendix.tex}

\bibliographystyle{abbrv}
\bibliography{LCSS2020}

\end{document}

%% file: preamble.tex
\usepackage[caption=false,font=footnotesize]{subfig}
\usepackage{url}
\usepackage[nocompress]{cite}

\usepackage{color}

\input{acronyms.tex}
\input{symbols.tex}

\usepackage{amsthm}
\theoremstyle{plain}
\ifx\c@thm\undefined
    \newtheorem{assumption}{Assumption}
\else
    \newtheorem{assumption}[thm]{Assumption}
\fi

\newenvironment{myproof}[1]{\noindent\hspace{2em}{\itshape Proof of #1: }}{\hspace*{\fill}~\QED\par\endtrivlist\unskip}

\allowdisplaybreaks[3]
\IEEEoverridecommandlockouts

%% file: acronyms.tex
\usepackage{acronym}
\acrodef{ALM}[ALM]{augmented Lagrangian method}
\acrodef{PDGM}[PDGM]{primal-dual gradient method}

%% file: symbols.tex
\usepackage{amsmath}
\usepackage{amssymb}
\usepackage{amsfonts}
\usepackage{mathtools}
\renewenvironment{align*}{\align}{\endalign}
\usepackage{autonum}

\newcommand{\restore@Environment}[1]{%
    \AtBeginDocument{%
        \csletcs{#1*}{#1}%
        \csletcs{end#1*}{end#1}%
    }%
}
\forcsvlist\restore@Environment{alignat,equation,gather,multline,flalign,align}

\usepackage{arydshln}

\newcommand{\norm}[1]{\left\Vert#1\right\Vert}
\DeclarePairedDelimiter{\cnorm}{\lVert}{\rVert}

\DeclareMathOperator*{\argmin}{arg\,min}

\newcommand{\deq}{\coloneqq}

\DeclareMathOperator*{\optmin}{\mathrm{min.}}
\DeclareMathOperator*{\optmax}{\mathrm{max.}}
\newcommand{\optst}{\mathrm{s.t.}}

\newcommand{\opt}{^{\star}}

\newcommand{\exact}{_{\mathrm{ex}}}
\newcommand{\extext}{``ex'' }
\newcommand{\gd}{_{\mathrm{gd}}}
\newcommand{\io}{_{\mathrm{io}}}

\newcommand{\sgvar}{s} 

\newcommand{\mup}{\mu_p}
\newcommand{\mug}{\mu_g}
\newcommand{\betapsi}{\beta_\phi}  
\newcommand{\betap}{\beta_p}
\newcommand{\betag}{\beta_g}
\newcommand{\betaf}{\beta_f}
\newcommand{\alphap}{\alpha_p}
\newcommand{\alphag}{\alpha_g}

\newcommand{\sect}{\mathsf{sec}} 
\newcommand{\Sin}{v} 
\newcommand{\Sout}{w} 
\newcommand{\Sstate}{\xi} 

\newcommand{\Bx}{B_\phi}  
\newcommand{\dxtoxi}{c_\phi}
\newcommand{\gainAE}{c_\xi}

%% file: affiliation.tex
\newcommand{\allthanks}{\thanks{The author is with the Department of Electrical and Computer Engineering, University of Illinois, Chicago, IL 60607. {\tt\scriptsize {hanshuo@uic.edu}}.}}

%% file: abstract.tex
We show that the primal-dual gradient method, also known as the gradient descent ascent method, for solving convex-concave minimax problems can be viewed as an inexact gradient method applied to the primal problem. The gradient, whose exact computation relies on solving the inner maximization problem, is computed approximately by another gradient method. To model the approximate computational routine implemented by iterative algorithms, we introduce the notion of \emph{dynamic inexact oracles}, which are discrete-time dynamical systems whose output asymptotically approaches the output of an exact oracle. We present a unified convergence analysis for dynamic inexact oracles realized by general first-order methods and demonstrate its use in creating new accelerated primal-dual algorithms. 

%% file: intro.tex
We consider algorithms for solving the unconstrained minimax problem
\begin{equation}
\optmin_{x\in\mathbb{R}^{n}}\quad\max_{y\in\mathbb{R}^{m}}L(x,y)\deq f(x)+y^{T}Ax-g(y).\label{eq:prob}
\end{equation}
We assume that $f$ is smooth and convex (but not necessarily strongly convex), $g$ is smooth and strongly convex, and $A\in\mathbb{R}^{m\times n}$ has full column rank. For convenience, we define $p(x)\deq\max_{y}L(x,y)$ and write problem~(\ref{eq:prob}) as
\begin{equation}
\optmin_{x}\quad p(x),\label{eq:prob_p}
\end{equation}
which we refer to as the \emph{primal problem}. We also define $d(y)\deq\min_{x}L(x,y)$ and refer to the problem
\begin{equation}
\optmax_{y}\quad d(y)\label{eq:prob_d}
\end{equation}
as the \emph{dual problem}. Under the given assumptions, it follows from standard results (see, e.g.,~\cite[Ch.~X]{hiriart-urruty_convex_1993}) in convex analysis that both $p$ and $-d$ are strictly convex (in fact, strongly convex). Therefore, the primal-dual optimal solution of problems~(\ref{eq:prob_p}) and~(\ref{eq:prob_d}) is unique, which we denote by $(x\opt,y\opt)$. 

The minimax problem~(\ref{eq:prob}) has a number of applications. For example, when $f(x)=-b^{T}x$ for some $b\in\mathbb{R}^{n}$, the dual problem~(\ref{eq:prob_d}) becomes equivalent to the equality-constrained convex optimization problem given by
\begin{equation}
\optmax_{y}\quad-g(y),\qquad\optst\quad A^{T}y=b.\label{eq:prob_eq_constr}
\end{equation}
Other applications include image processing~\cite{chambolle_first-order_2011} and empirical risk minimization~\cite{zhang_stochastic_2017}. More broadly, when the function $L$ is a general convex-concave function, the minimax problem formulation also arises in game theory~\cite{myerson_game_2013} and robust optimization~\cite{ben-tal_robust_2009}.

One important algorithm for computing the primal-dual optimal solution $(x\opt,y\opt)$ is the \emph{\ac{PDGM}}:
\begin{equation}
\begin{aligned}x^{k+1} & =x^{k}-\eta_{1}\nabla_{1}L(x^{k},y^{k})\\
y^{k+1} & =y^{k}+\eta_{2}\nabla_{2}L(x^{k},y^{k}),
\end{aligned}
\label{eq:pdgm_original}
\end{equation}
where $\eta_{1}$ and $\eta_{2}$ are step sizes, and $\nabla_{1}L(x^{k},y^{k})=\nabla f(x^{k})+A^{T}y^{k}$ and $\nabla_{2}L(x^{k},y^{k})=Ax^{k}-\nabla g(y^{k})$ are the partial derivatives of $L$ with respect to the first and second arguments, respectively. The \ac{PDGM} is also known by various other names such as the Arrow\textendash Hurwicz gradient method~\cite[p.~155]{arrow_studies_1958} and the (simultaneous) gradient descent ascent method (see, e.g.,~\cite{daskalakis_limit_2018}). It has also been generalized to the case where $L$ is non-differentiable~\cite{nedic_subgradient_2009} and the case where the dynamics in~(\ref{eq:pdgm_original}) are in continuous time~\cite{feijer_stability_2010,cherukuri_saddle-point_2017,qu_exponential_2019}. Convergence of the \ac{PDGM} has been studied extensively in the literature. Under the assumption we made on $f$, $g$, and $A$, it has been shown~\cite{du_linear_2019} that the \ac{PDGM} converges exponentially to the optimal solution $(x\opt,y\opt$).

Because the update rule~(\ref{eq:pdgm_original}) of the \ac{PDGM} performs gradient descent/ascent on the primal/dual variable, a natural question arises as to whether these gradient updates can be substituted by other first-order methods (e.g., Nesterov's accelerated gradient method) to create new primal-dual algorithms. Our paper attempts to address this question by providing a unified convergence analysis that allows the gradient updates to be replaced by a class of first-order methods. The analysis hinges on an alternative view of the \ac{PDGM}: We show that the \ac{PDGM} is equivalent to applying an inexact gradient method to the primal problem~(\ref{eq:prob_p}), where the gradient $\nabla p$ is computed approximately by a \emph{dynamic inexact oracle} (see Definition~\ref{def:DIO}). A dynamic inexact oracle is only required to compute the exact gradient asymptotically, and the transient results of such an oracle may be inexact. For the case of the \ac{PDGM}, the inexact oracle is realized by running one iteration of gradient descent with warm starts (see Section~\ref{subsec:pdgm_as_inexact}). This abstract view using dynamic inexact oracles leads to a unified convergence analysis that does not rely on the detailed realization of the oracle. %

\paragraph*{Contribution}

While the notion of inexact oracles has long existed in the study of optimization algorithms, including approximating the gradient mapping (see, e.g., \cite[Ch.~3.3]{bertsekas_convex_2015}) and the proximal operator~\cite{rockafellar_monotone_1976}, these inexact oracles are static mappings and hence less general than our proposed notion of dynamic inexact oracles, which are permitted to have internal states and are necessary for modeling warm starts used in iterative algorithms. The introduction of dynamics also demands a new analysis for understanding the dynamical interaction between the gradient method and the inexact oracle used therein. By modeling the dynamical interaction as a feedback interconnection of two dynamical systems, we derive a convergence analysis (Theorems~\ref{thm:pdgm_convergence} and~\ref{thm:io_convergence}) using the small-gain principle. The convergence analysis also enables us to build new primal-dual algorithms by simply changing the realization of the inexact oracle used in \ac{PDGM} to other first-order methods in a ``plug-and-play'' manner.

%% file: prelim.tex
For a vector $x$, we denote by $\norm{x}$ its $\ell_{2}$-norm and $\norm{x}_{P}\deq(x^{T}Px)^{1/2}$ its $P$-quadratic norm, where $P$ is a positive definite matrix (written as $P\succ0$). For a function $f(\cdot,\cdot)$ with two arguments, we denote by $\nabla_{i}f$ ($i=1,2$) the partial derivative of $f$ with respect to the $i$th argument. Unless noted otherwise, we reserve the use of superscripts for indexing an infinite sequence $\{x^{k}\}_{k=0}^{\infty}$. 

For a real-valued function $f$, we denote by $f^{*}$ its convex conjugate, defined by $f^{*}(s)\deq\sup_{x}\{s^{T}x-f(x)\}$. We denote by $\mathcal{S}(\mu,\beta)$ the set of $\mu$-strongly convex and $\beta$-smooth functions. By convention, we use $\mathcal{S}(0,\beta)$ to denote the set of $\beta$-smooth and convex functions. Recall the following basic properties of functions in $\mathcal{S}(\mu,\beta)$. 
\begin{prop}
[Basic properties]\label{prop:smooth_and_scvx}If $f\in\mathcal{S}(\mu,\beta)$, then
\end{prop}
\begin{enumerate}
\item $(x-y,\nabla f(x)-\nabla f(y))\in\sect(\mu,\beta)$ for all $x$ and $y$, where 
\begin{multline*}
\sect(\mu,\beta)\deq\left\{ (\Sin,\Sout)\colon\left[\begin{array}{c}
\Sin\\
\Sout
\end{array}\right]^{T}\right.\\
\left.\left[\begin{array}{cc}
-2\mu\beta I & (\mu+\beta)I\\
(\mu+\beta)I & -2I
\end{array}\right]\left[\begin{array}{c}
\Sin\\
\Sout
\end{array}\right]\right\} \geq0
\end{multline*}
is called the \emph{sector constraint.}
\end{enumerate}
Furthermore, if $\mu>0$, then
\begin{enumerate}
\item [2)]\setcounter{enumi}{2}$f^{*}\in\mathcal{S}(1/\beta,1/\mu)$; 
\item $\nabla f$ is invertible and $(\nabla f)^{-1}=\nabla f^{*}$, where $\nabla f^{*}$ is the gradient of $f^{*}$;
\end{enumerate}
A proof of item 1 can be found in~\cite[Thm.~2.1.12]{nesterov_introductory_2004}. Proofs of items 2 and 3 can be found in~\cite[Ch.~X]{hiriart-urruty_convex_1993}.

%% file: oracle.tex
We begin by considering another way to solve the primal problem~(\ref{eq:prob_p}) by directly applying the gradient method%
. By allowing inexact gradient computation, we reveal that the \ac{PDGM} can be viewed alternatively as an inexact gradient method applied to the primal problem%
. An abstraction of the inexact gradient computation leads to the definition of dynamic inexact oracles, the central topic of study in this paper. 

\subsection{The \ac{PDGM} as inexact gradient descent \label{subsec:pdgm_as_inexact}}

Consider solving the primal problem~(\ref{eq:prob_p}) using the gradient method:
\begin{equation}
x^{k+1}\exact=x^{k}\exact-\eta_{1}\nabla p(x^{k}\exact),\label{eq:primal_grad_original}
\end{equation}
where $\eta_{1}$ is the step size. (The subscript \extext stands for \emph{exact}, in comparison to the inexact gradient method to be presented shortly.) Define $\tilde{g}(y,x)\deq f(x)-L(x,y)=g(y)-y^{T}Ax$. Using Danskin's theorem (see, e.g.,~\cite[p.~245]{bertsekas_convex_2003}), we obtain $\nabla p(x^{k}\exact)=\nabla f(x^{k}\exact)+A^{T}y^{k}\exact$, where $y^{k}\exact=\argmin_{y}\tilde{g}(y,x^{k}\exact)$ (unique because $\tilde{g}$ is strongly convex). Therefore, the gradient method~(\ref{eq:primal_grad_original}) can be rewritten as\begin{subequations}\label{eq:primal_grad}
\begin{align}
y^{k}\exact & =\argmin_{y}\tilde{g}(y,x^{k}\exact)\label{eq:primal_grad_y}\\
x^{k+1}\exact & =x^{k}\exact-\eta_{1}(\nabla f(x^{k}\exact)+A^{T}y^{k}\exact).\label{eq:primal_grad_x}
\end{align}
\end{subequations}
\begin{rem}
The equality-constrained optimization problem~(\ref{eq:prob_eq_constr}) can be viewed as the dual problem of~(\ref{eq:prob}) for $L(x,y)=-b^{T}x+y^{T}Ax-(g(y)+\frac{\eta_{1}}{2}\cnorm{A^{T}y-b}^{2})$. In this case, the functions $f$ and $\tilde{g}$ are given by $f(x)=-b^{T}x$ and $\tilde{g}(y,x)=g(y)+\frac{\eta_{1}}{2}\cnorm{A^{T}y-b}^{2}-y^{T}Ax$, and the recursion (\ref{eq:primal_grad}) becomes
\begin{align*}
y^{k}\exact & =\argmin_{y}\left\{ g(y)-y^{T}Ax^{k}\exact+\frac{\eta_{1}}{2}\norm{A^{T}y-b}^{2}\right\} \\
x^{k+1}\exact & =x^{k}\exact-\eta_{1}(A^{T}y^{k}\exact-b),
\end{align*}
which recovers the augmented Lagrangian method (see, e.g.,~\cite[p.~262]{bertsekas_convex_2015}).
\end{rem}
Under appropriate choice of the step size $\eta_{1}$, the sequence $\{x^{k}\exact\}$ generated by~(\ref{eq:primal_grad}) converges to the optimal solution $x\opt$ of the primal problem~(\ref{eq:prob_p}). However, because the gradient mapping $\nabla p$ depends on $y\exact$, each iteration requires solving the minimization problem in~(\ref{eq:primal_grad_y}). This is undesired because the minimization problem does not generally admit a closed-form solution. 

We now show that the \ac{PDGM} can be derived from~(\ref{eq:primal_grad}) by allowing the minimization problem in~(\ref{eq:primal_grad_y}) to be solved approximately. Suppose the approximate solution, denoted by $\{y^{k}\}$, is generated by applying one iteration of the gradient method (with step size $\eta_{2}$) to the problem in~(\ref{eq:primal_grad_y}). This yields 
\begin{equation}
y^{k+1}=y^{k}-\eta_{2}\nabla_{1}\tilde{g}(y^{k},x^{k})=y^{k}+\eta_{2}(Ax^{k}-\nabla g(y^{k})).\label{eq:y_approx}
\end{equation}
Note that the update rule~(\ref{eq:y_approx}) uses a warm start: It uses the approximate solution $y^{k}$ at iteration $k$ to initialize iteration $k+1$. The approximate solution $\{y^{k}\}$ is then used in place of $\{y^{k}\exact\}$ in~(\ref{eq:primal_grad_x}), yielding
\begin{equation}
x^{k+1}=x^{k}-\eta_{1}(\nabla f(x^{k})+A^{T}y^{k}),\label{eq:x_approx}
\end{equation}
where we have replaced $x\exact$ with $x$ to distinguish from the sequence generated by the exact gradient method. It can be seen that the update rules~(\ref{eq:y_approx}) and~(\ref{eq:x_approx}) recover the \ac{PDGM}; namely, the \ac{PDGM} can be viewed as an inexact gradient method applied to the primal problem. 

\subsection{Dynamic inexact oracles\label{subsec:dio_def}}

It is not difficult to imagine that the gradient method~(\ref{eq:y_approx}) is not the only iterative algorithm for generating an approximate solution to the minimization problem in~(\ref{eq:primal_grad_y}), which is needed for computing the gradient $\nabla p$. To facilitate discussion, we introduce the notion of dynamic inexact oracles as a high-level description of iterative algorithms used for approximation. 
\begin{defn}
[Dynamic inexact oracles] \label{def:DIO}A (discrete-time) dynamical system $\mathcal{G}$ is called a \emph{dynamic inexact oracle} for computing a mapping $\phi$ if for any input sequence $u=\{u^{k}\}_{k=0}^{\infty}$ converging to $u\opt$, the output $\mathcal{G}u$ converges to $\phi(u\opt)$. 
\end{defn}
If $\mathcal{G}$ is a dynamic inexact oracle, even when the input sequence $u\equiv u\opt$ is constant, the output of $\mathcal{G}$ is not required to immediately match the exact oracle output $\phi(u\opt)$, hence the term \emph{inexact}; the only requirement is that $\mathcal{G}$ must compute $\phi(u\opt)$ asymptotically. 

In the remainder of this paper, we focus on dynamic inexact oracles that approximately solve the optimization problem in~(\ref{eq:primal_grad_y}). Denote by $\{x^{k}\}$ and $\{y^{k}\}$ the input and output of the oracle, respectively. For any $\{x^{k}\}$ converging to $x\opt$, the output of the inexact oracle must asymptotically converge to the optimal solution $y\opt=\argmin_{y}\tilde{g}(y,x\opt)$. We will show in Section~\ref{subsec:oracle_gd} that the update rule~(\ref{eq:y_approx}) given by the gradient method is one such inexact oracle. Furthermore, by constructing the inexact oracle from different first-order optimization algorithms, we can create new primal-dual first-order methods beyond the \ac{PDGM} (see Section~\ref{subsec:oracle_1st}). 

The notion of dynamic inexact oracles is fundamentally different from the inexact oracles studied in the existing literature, which are static inexact oracles. For a static oracle, the output of the oracle at any iteration $k$ only depends on the instantaneous input $u^{k}$. Incorporating dynamics into inexact oracles is necessary because a static oracle is not able to model iterative optimization algorithms with warm starts, in which the solution during the current iteration needs to be memorized to initialize the next iteration such as in~(\ref{eq:y_approx}). One example of static inexact oracles is approximate gradient mappings used in first-order methods, such as in the $\epsilon$-(sub)gradient method (see, e.g.,~\cite[Ch.~3.3]{bertsekas_convex_2015}). Other examples include approximate proximal operators used in the proximal point algorithm~\cite[p.~880]{rockafellar_monotone_1976} and in the Douglas\textendash Rachford splitting method \cite[Thm.~8]{eckstein_douglasrachford_1992}. A general treatment of static inexact oracles in first-order methods can be found in~\cite{devolder_first-order_2014}.

%% file: analysis.tex
We show that the convergence of gradient methods with dynamic inexact oracles can be analyzed by viewing it as a feedback interconnection of two dynamical systems. By applying the small-gain principle, we present a unified convergence analysis that only depends on the input-output behavior of the inexact oracle. We begin with the oracle realized by gradient descent, after which we extend the analysis to oracles realized by general first-order methods.

We shall make the following assumptions on $f$, $g$, and $A$:\begin{assumption}\label{assu:fgA}

Let $f$, $g$, and $A$ in the minimax problem~(\ref{eq:prob}) be such that $f\in\mathcal{S}(0,\betaf)$, $g\in\mathcal{S}(\mug,\betag)$, and $A$ has full column rank. \end{assumption}

Let $\sigma_{\max}$ and $\sigma_{\min}$ be the maximum and minimum singular values of $A$, respectively. Recall that the primal objective function $p$ is given by $p(x)=\max_{y}L(x,y)=f(x)+g^{*}(Ax)$. From Proposition~\ref{prop:smooth_and_scvx}, we have $p\in\mathcal{S}(\mup,\betap)$, where $\mup=\sigma_{\min}^{2}/\betag$ and $\betap=\sigma_{\max}^{2}/\mug+\betaf$. 

\subsection{The oracle based on gradient descent\label{subsec:oracle_gd}}

We begin by showing that the recursion~(\ref{eq:y_approx}) based on gradient descent, which can be viewed as a dynamical system $\mathcal{G}\gd$ with input $x$ and output $y$, is indeed an inexact oracle for computing the optimal solution of the minimization problem in~(\ref{eq:primal_grad_y}). The optimality condition of the minimization problem gives
\[
0=\nabla_{1}\tilde{g}(y^{k}\exact,x^{k}\exact)=\nabla g(y^{k}\exact)-Ax^{k}\exact.
\]
Because $g\in\mathcal{S}(\mug,\betag)$, using Proposition~\ref{prop:smooth_and_scvx}, we obtain 
\begin{equation}
y^{k}\exact=(\nabla g)^{-1}(Ax^{k}\exact)=\nabla g^{*}(Ax^{k}\exact)\eqqcolon\phi(x^{k}\exact)\label{eq:def_phi}
\end{equation}
In other words, we need to show that $\mathcal{G}\gd$ asymptotically computes the mapping $\phi$. 
\begin{prop}
\label{prop:G_gd}Let $\mathcal{G}\gd$ be a dynamical system defined by~(\ref{eq:y_approx}) and $\eta_{2}\in(0,2/(\mug+\betag)]$. Suppose the input $\{x^{k}\}$ converges to $x\opt$. Then the output $\{y^{k}\}$ of $\mathcal{G}\gd$ converges to $\phi(x\opt)=\nabla g^{*}(Ax\opt)$. 
\end{prop}
To prove Proposition~\ref{prop:G_gd}, we need to make use of the following lemma, which plays an important role in establishing the convergence of first-order methods. See Appendix for a proof, which is a straightforward consequence of standard results in convex optimization.
\begin{lem}
\label{lem:sec_contr}Let $\mu$ and $\beta$ be constants satisfying $0<\mu\leq\beta$, and $\alpha=\mu\beta/(\mu+\beta)$. Suppose $(\Sstate,\Sout)\in\sect(\mu,\beta)$. Then for any $\eta\in(0,2/(\mu+\beta)]$, we have $\cnorm{\Sstate-\eta\Sout}\leq\rho\cnorm{\Sstate}$, where $\rho=1-\alpha\eta\in[0,1)$. 
\end{lem}
\begin{myproof}{Proposition~\ref{prop:G_gd}}Define $y\opt\deq\nabla g^{*}(Ax\opt)$, and rewrite~(\ref{eq:y_approx}) as 
\[
y^{k+1}-y\opt=(y^{k}-y\opt)-\eta_{2}(\nabla g(y^{k})-Ax\opt)+\eta_{2}A(x^{k}-x\opt).
\]
Because $g\in\mathcal{S}(\mug,\betag)$ and $Ax\opt=\nabla g(\nabla g^{*}(Ax\opt))=\nabla g(y\opt)$, we have $(y^{k}-y\opt,\nabla g(y^{k})-Ax\opt)\in\sect(\mug,\betag)$. By Lemma~\ref{lem:sec_contr}, there exists $\rho\in[0,1)$ such that
\[
\cnorm{y^{k+1}-y\opt}\leq\rho\cnorm{y^{k}-y\opt}+\eta_{2}\cnorm{A(x^{k}-x\opt)}.
\]
The result then follows as a consequence of input-to-state stability~\cite[p.~192]{khalil_nonlinear_2002}.\end{myproof}

\begin{figure}
\begin{centering}
\includegraphics[width=0.75\columnwidth]{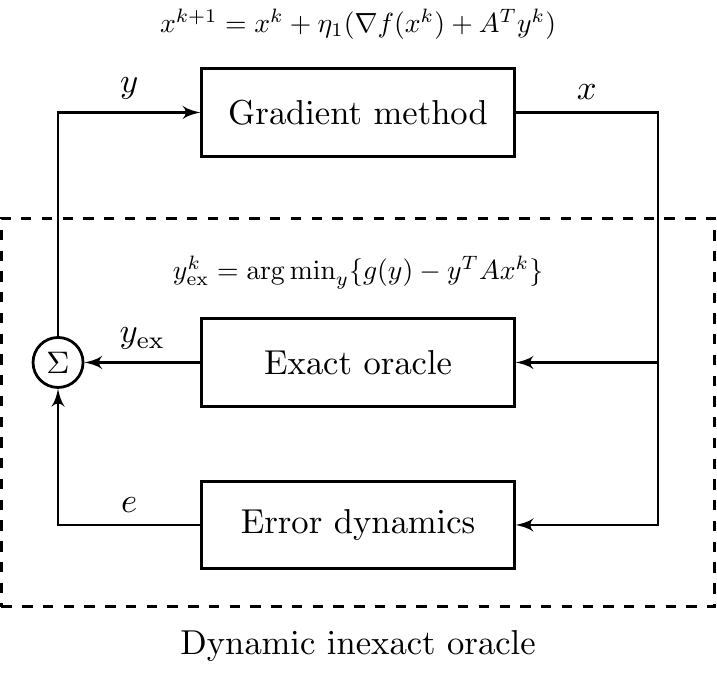}
\par\end{centering}
\caption{The gradient method with a dynamic inexact oracle. The difference between the inexact and the exact oracles is characterized by the additive error dynamics. \label{fig:gd_dio}}

\end{figure}

Let us now analyze the convergence of the gradient method~(\ref{eq:x_approx}) with the dynamic inexact oracle $\mathcal{G}\gd$ defined by~(\ref{eq:y_approx}). For convenience, we define the error $e$ such that $e^{k}\deq y^{k}-\nabla g^{*}(Ax^{k})$ and rewrite~(\ref{eq:y_approx}) and~(\ref{eq:x_approx}) as\begin{subequations}\label{eq:gd_recursion_xe}
\begin{align}
x^{k+1} & =x^{k}-\eta_{1}\nabla p(x^{k})-\eta_{1}A^{T}e^{k}.\label{eq:x_approx_e}\\
e^{k+1} & =e^{k}-\eta_{2}\nabla_{1}\tilde{g}(y^{k},x^{k})-[\nabla g^{*}(Ax^{k+1})-\nabla g^{*}(Ax^{k})].\label{eq:e_approx}
\end{align}
\end{subequations} Although the recursion~(\ref{eq:x_approx_e}) converges when the error $e\equiv0$, and the recursion~(\ref{eq:e_approx}) converges when $x\equiv x\opt$ (Proposition~\ref{prop:G_gd}), the joint recursion~(\ref{eq:gd_recursion_xe}) is not guaranteed to converge. Indeed, the joint recursion~(\ref{eq:gd_recursion_xe}) can be viewed as a feedback interconnection of two dynamical systems~(\ref{eq:x_approx_e}) and~(\ref{eq:e_approx}) as illustrated in Fig.~\ref{fig:gd_dio}, and it is well known in control theory that a feedback connection of two internally stable systems may be unstable.

A powerful method for analyzing the stability of feedback interconnections of dynamical systems is the small-gain principle. The small-gain principle can take various forms depending on the specific setup. The following is what we will use in this paper. See Appendix for a detailed proof.
\begin{lem}
[Small-gain] \label{lem:small-gain}Let $\{\sgvar_{1}^{k}\}$ and $\{\sgvar_{2}^{k}\}$ be two nonnegative real-valued sequences satisfying
\begin{align*}
\sgvar_{1}^{k+1} & \leq\gamma_{11}\sgvar_{1}^{k}+\gamma_{12}\sgvar_{2}^{k}\\
\sgvar_{2}^{k+1} & \leq\gamma_{21}\sgvar_{1}^{k}+\gamma_{22}\sgvar_{2}^{k}
\end{align*}
for some nonnegative constants $\gamma_{ij}$ $(i,j=1,2)$. Then, both $\{\sgvar_{1}^{k}\}$ and $\{\sgvar_{2}^{k}\}$ converge exponentially to $0$ if $\gamma_{11}<1$, $\gamma_{22}<1$, and $\gamma_{12}\gamma_{21}<(1-\gamma_{11})(1-\gamma_{22})$. 
\end{lem}

The small-gain lemma (Lemma~\ref{lem:small-gain}) shows that, in order for the feedback interconnection of two (nonnegative) systems to be stable, aside from the stability of individual systems ($\gamma_{11}<1$ and $\gamma_{22}<1$), the coupling coefficients $\gamma_{12}$ and $\gamma_{21}$ must be small enough. We now apply the small-gain lemma to establish the convergence of~(\ref{eq:gd_recursion_xe}). 
\begin{thm}
\label{thm:pdgm_convergence}Consider the gradient method given by~(\ref{eq:x_approx}), where $\{y^{k}\}$ is given by the dynamic inexact oracle $\mathcal{G}\gd$ defined by~(\ref{eq:y_approx}) with $\eta_{2}\in(0,2/(\mug+\betag)]$. Suppose $f$, $g$, and $A$ satisfy Assumption~\ref{assu:fgA}, and let $\betapsi=\sigma_{\max}/\mug$, $\alphap=\mup\betap/(\mup+\betap)$, and $\alphag=\mug\betag/(\mug+\betag)$. Then, for any $\eta_{1}$ satisfying
\begin{equation}
0<\eta_{1}<\min\left\{ \frac{\alphap\alphag\eta_{2}}{\sigma_{\max}\betapsi(\alphap+\betap)},\frac{2}{\mup+\betap}\right\} ,\label{eq:eta1_gd}
\end{equation}
the sequences $\{x^{k}\}$ and $\{y^{k}\}$ converge exponentially to the primal and dual optimal solutions $x\opt$ and $y\opt$, respectively. 
\end{thm}
\begin{IEEEproof}
Denote by $e\opt$ the steady-state value of $e$. Then, we have $e\opt=y\opt-\nabla g^{*}(Ax\opt)=0$. Define $\hat{x}^{k}\deq x^{k}-x\opt$ and $\hat{e}^{k}\deq e^{k}-e\opt$. We can rewrite~(\ref{eq:gd_recursion_xe}) as 
\begin{align*}
\hat{x}^{k+1} & =\hat{x}^{k}-\eta_{1}\nabla p(x^{k})-\eta_{1}A^{T}\hat{e}^{k}\\
\hat{e}^{k+1} & =\hat{e}^{k}-\eta_{2}\nabla_{1}\tilde{g}(y^{k},x^{k})-[\nabla g^{*}(Ax^{k+1})-\nabla g^{*}(Ax^{k})].
\end{align*}
Because $p\in\mathcal{S}(\mup,\betap)$ and $g\in\mathcal{S}(\mug,\betag)$, from Proposition~\ref{prop:smooth_and_scvx}, we have $(\hat{x}^{k},\nabla p(x^{k}))=(x^{k}-x\opt,\nabla p(x^{k})-\nabla p(x\opt))\in\sect(\mup,\betap)$ and $(\hat{e}^{k},\nabla_{1}\tilde{g}(y^{k},x^{k}))=(y^{k}-\nabla g^{*}(Ax^{k}),\nabla g(y^{k})-Ax^{k})\in\sect(\mug,\betag),$ where we have used the fact $Ax^{k}=\nabla g(\nabla g^{*}(Ax^{k}))$. Applying Lemma~\ref{lem:sec_contr}, since $0<\eta_{1}\leq2/(\mup+\betap)$, we have 
\begin{align}
\cnorm{\hat{x}^{k+1}} & \leq\cnorm{\hat{x}^{k}-\eta_{1}\nabla p(x^{k})}+\norm{\eta_{1}A^{T}\hat{e}^{k}}\nonumber \\
 & \leq\rho_{1}\cnorm{\hat{x}^{k}}+\eta_{1}\sigma_{\max}\cnorm{\hat{e}^{k}},\label{eq:gd_x_ineq}
\end{align}
where $\rho_{1}=1-\alphap\eta_{1}$; similarly, we also obtain 
\begin{align}
\cnorm{\hat{e}^{k+1}} & \leq\rho_{2}\cnorm{\hat{e}^{k}}+\betapsi\cnorm{x^{k+1}-x^{k}}\nonumber \\
 & =\rho_{2}\cnorm{\hat{e}^{k}}+\betapsi\cnorm{-\eta_{1}\nabla p(x^{k})-\eta_{1}A^{T}\hat{e}^{k}}\nonumber \\
 & \leq\eta_{1}\betapsi\betap\cnorm{\hat{x}^{k}}+(\rho_{2}+\eta_{1}\betapsi\sigma_{\max})\cnorm{\hat{e}^{k}},\label{eq:gd_e_ineq}
\end{align}
where $\rho_{2}\in1-\alphag\eta_{2}$, and $\betapsi=\sigma_{\max}/\mu_{g}$ is the Lipschitz constant of the mapping $\phi\colon x^{k}\mapsto\nabla g^{*}(Ax^{k})$. The relationship given by~(\ref{eq:gd_x_ineq}) and~(\ref{eq:gd_e_ineq}) allows us to apply the small-gain lemma (Lemma~\ref{lem:small-gain}) and derive the condition~(\ref{eq:eta1_gd}) for both $\hat{x}$ and $\hat{e}$ to converge exponentially to $0$, i.e., $x^{k}\to x\opt$ and $y^{k}\to\nabla g^{*}(Ax\opt)=y\opt$ as required. 
\end{IEEEproof}

Although exponential convergence of the \ac{PDGM} has already been established~\cite{du_linear_2019}, the technique used in the proof of Theorem~\ref{thm:pdgm_convergence} is different from what is used in the existing literature. The proof reveals two attractive features of the small-gain principle in the analysis of the inexact gradient method. First, it is capable of incorporating existing convergence results, i.e., internal stability of the gradient dynamics~(\ref{eq:x_approx_e}) and~(\ref{eq:e_approx}) as manifested in Lemma~\ref{lem:sec_contr}. This avoids the need of finding a Lyapunov function from scratch; in comparison, typical convergence proofs of first-order algorithms in the literature involve constructing a Lyapunov function, which is often nontrivial except for the simplest algorithms. Second, the small-gain analysis only relies on a coarse description of the input-output behavior such as what is given in~(\ref{eq:gd_e_ineq}). Therefore, when the dynamic inexact oracle is realized by another iterative algorithm $\mathcal{G}\io$, the small-gain analysis can be readily applied as long as a relationship between the input $x$ and the error $e$ of $\mathcal{G}\io$ similar to~(\ref{eq:gd_e_ineq}) can be derived (which, incidentally, often makes use of the fact that $\mathcal{G}\io$ is a dynamic inexact oracle and hence internally stable). The ``plug-and-play'' nature of this approach allows us to easily generalize the analysis to a wide range of dynamic inexact oracles, which we will discuss shortly in Section~\ref{subsec:oracle_1st}. 

\subsection{Oracles based on general first-order algorithms\label{subsec:oracle_1st}}

As we pointed out in Section~\ref{subsec:dio_def}, dynamical inexact oracles for solving the minimization problem in~(\ref{eq:primal_grad_y}) can be constructed from iterative optimization algorithms. Inspired by the work in~\cite{van_scoy_fastest_2018}, we consider inexact oracles constructed from algorithms in the following state-space form:
\begin{equation}
\begin{aligned}\xi^{k+1} & =A\io\xi^{k}+B\io\nabla F(v^{k})\\
v^{k} & =C\io\xi^{k},\qquad z^{k}=E\io\xi^{k},
\end{aligned}
\label{eq:opt_alg_compact}
\end{equation}
where $A\io$, $B\io$, $C\io$, and $E\io$ are given by 
\begin{align*}
A\io & =\left[\begin{array}{cc}
(1+c_{1})I & -c_{1}I\\
I & 0
\end{array}\right],\quad B\io=\left[\begin{array}{c}
-\eta_{2}I\\
0
\end{array}\right],\\
C\io & =\begin{bmatrix}(1+c_{2})I & -c_{2}I\end{bmatrix},\quad E\io=\begin{bmatrix}(1+c_{3})I & -c_{3}I\end{bmatrix}.
\end{align*}
Here, $F$ is the (convex) objective function to be minimized, $\xi=(\xi_{1},\xi_{2})$ is the state, $v$ is the feedback output, $z$ is the output of the algorithm, $\eta_{2}$ is the step size, and $c_{1}$, $c_{2}$, and $c_{3}$ are constants. The form~(\ref{eq:opt_alg_compact}) captures a number of important first-order optimization algorithms. For example, setting $c_{1}=c_{2}=c_{3}=0$ recovers the gradient method, and setting $c_{1}=c_{2}\neq0$ and $c_{3}=0$ recovers Nesterov's accelerated gradient method. Interested readers can refer to~\cite[Table I]{van_scoy_fastest_2018} for more examples. %

Similar to $\mathcal{G}\gd$ defined by~(\ref{eq:y_approx}), we construct a dynamic inexact oracle $\mathcal{G}\io$ by replacing $\nabla F(v^{k})$ and output $z^{k}$ in~(\ref{eq:opt_alg_compact}) with $\nabla_{1}\tilde{g}(v^{k},x^{k})=\nabla g(v^{k})-Ax^{k}$ and $y^{k}$, respectively:
\begin{equation}
\begin{aligned}\xi^{k+1} & =A\io\xi^{k}+B\io\left[\nabla g(v^{k})-Ax^{k}\right]\\
v^{k} & =C\io\xi^{k},\qquad y^{k}=E\io\xi^{k},
\end{aligned}
\label{eq:dio_1st}
\end{equation}
where $\{x^{k}\}$ and $\{y^{k}\}$ are the input and output of $\mathcal{G}\io$. 

Similar to Lemma~\ref{lem:sec_contr}, we shall make the following assumption on the algorithm given in~(\ref{eq:opt_alg_compact}).\begin{assumption}

\label{assm:ABC_io}Let $\mu$ and $\beta$ be constants satisfying $0<\mu\leq\beta$. Then there exist $P\succ0$, $\eta_{2}>0$, and $\rho_{2}\in[0,1)$ such that 
\[
\cnorm{A\io\Sstate+B\io\Sout}_{P}\leq\rho_{2}\cnorm{\Sstate}_{P}
\]
for all $\Sout$ satisfying $(\Sin,\Sout)\in\sect(\mu,\beta)$, where $\Sin=C\io\Sstate$. \end{assumption}

Assumption~\ref{assm:ABC_io} ensures that $\mathcal{G}\io$ is a dynamic inexact oracle that asymptotically computes the mapping $\phi$ defined in~(\ref{eq:def_phi}). The proof is similar to that of Proposition~\ref{prop:G_gd}. Recall that we can recover the gradient method by setting $c_{1}=c_{2}=c_{3}=0$ in~(\ref{eq:opt_alg_compact}). In this case, the second component $\xi_{2}$ of $\xi$ becomes irrelevant and can be dropped, so that we obtain (with an abuse of notion) $A\io=I$, $B\io=-\eta_{2}I$, and $C\io=I$. Therefore, by Lemma~\ref{lem:sec_contr}, the gradient method satisfies Assumption~\ref{assm:ABC_io} with $P=I$. For other first-order algorithms, while we we are unable to provide conditions under which Assumption~\ref{assm:ABC_io} holds, numerical methods~\cite[Figs.~3 and 5]{lessard_analysis_2016} have been used to show the existence of $P$, $\eta_{2}$, and $\rho$ for both Nesterov's accelerated gradient method and the heavy-ball method, at least when $\beta/\mu$ is small.

Convergence of the inexact gradient method~(\ref{eq:x_approx}) using a dynamic inexact oracle $\mathcal{G}\io$ can be established using a small-gain analysis similar to the proof of Theorem~\ref{thm:pdgm_convergence}. Details of the proof can be found in Appendix. 
\begin{thm}
\label{thm:io_convergence}Consider the gradient method given by~(\ref{eq:x_approx}), where $\{y^{k}\}$ is given by a dynamic inexact oracle $\mathcal{G}\io$ of the form~(\ref{eq:dio_1st}). Suppose $f$, $g$, and $A$ satisfy Assumption~\ref{assu:fgA}, and $A\io$, $B\io$, and $C\io$ satisfy Assumption~\ref{assm:ABC_io}. Then there exists $\eta_{1}$ such that $\{x^{k}\}$ and $\{y^{k}\}$ converge exponentially to the primal and dual optimal solutions $x\opt$ and $y\opt$, respectively. 
\end{thm}
As an application of Theorem~\ref{thm:io_convergence}, we give a convergence result for the case where $\mathcal{G}\io$ is realized by Nesterov's accelerated gradient method.
\begin{cor}
\label{cor:pd-nesterov}Let $\gamma=(\sqrt{\betag}-\sqrt{\mug})/(\sqrt{\betag}+\sqrt{\mug})$ and $\eta_{2}=1/\betag$. Consider the gradient method given by~(\ref{eq:x_approx}), where $\{y^{k}\}$ is given by a dynamic inexact oracle realized by Nesterov's accelerated gradient method:
\begin{equation}
\begin{aligned}y^{k+1} & =v^{k}-\eta_{2}(\nabla g(v^{k})-Ax^{k})\\
v^{k+1} & =(1+\gamma)y^{k+1}-\gamma y^{k}.
\end{aligned}
\label{eq:dio_nesterov}
\end{equation}
Suppose $f$, $g$, and $A$ satisfy Assumption~\ref{assu:fgA}. Then there exists $\eta_{1}$ such that $\{x^{k}\}$ and $\{y^{k}\}$ converge exponentially to the primal and dual optimal solutions $x\opt$ and $y\opt$, respectively, when $\betag/\mug$ is small enough.
\end{cor}
\begin{IEEEproof}
The recursion~(\ref{eq:dio_nesterov}) can be derived from~(\ref{eq:dio_1st}) by setting $c_{1}=c_{2}=\gamma$ and $c_{3}=0$ followed by eliminating $\xi$. Under the given choice of $\gamma$ and $\eta_{2}$, it has been shown in~\cite[Fig.~3]{lessard_analysis_2016} that Assumption~\ref{assm:ABC_io} holds when $\betag/\mug$ is small enough. The corollary then follows from Theorem~\ref{thm:io_convergence}.
\end{IEEEproof}

\begin{figure}
\begin{centering}
\includegraphics[width=1\columnwidth]{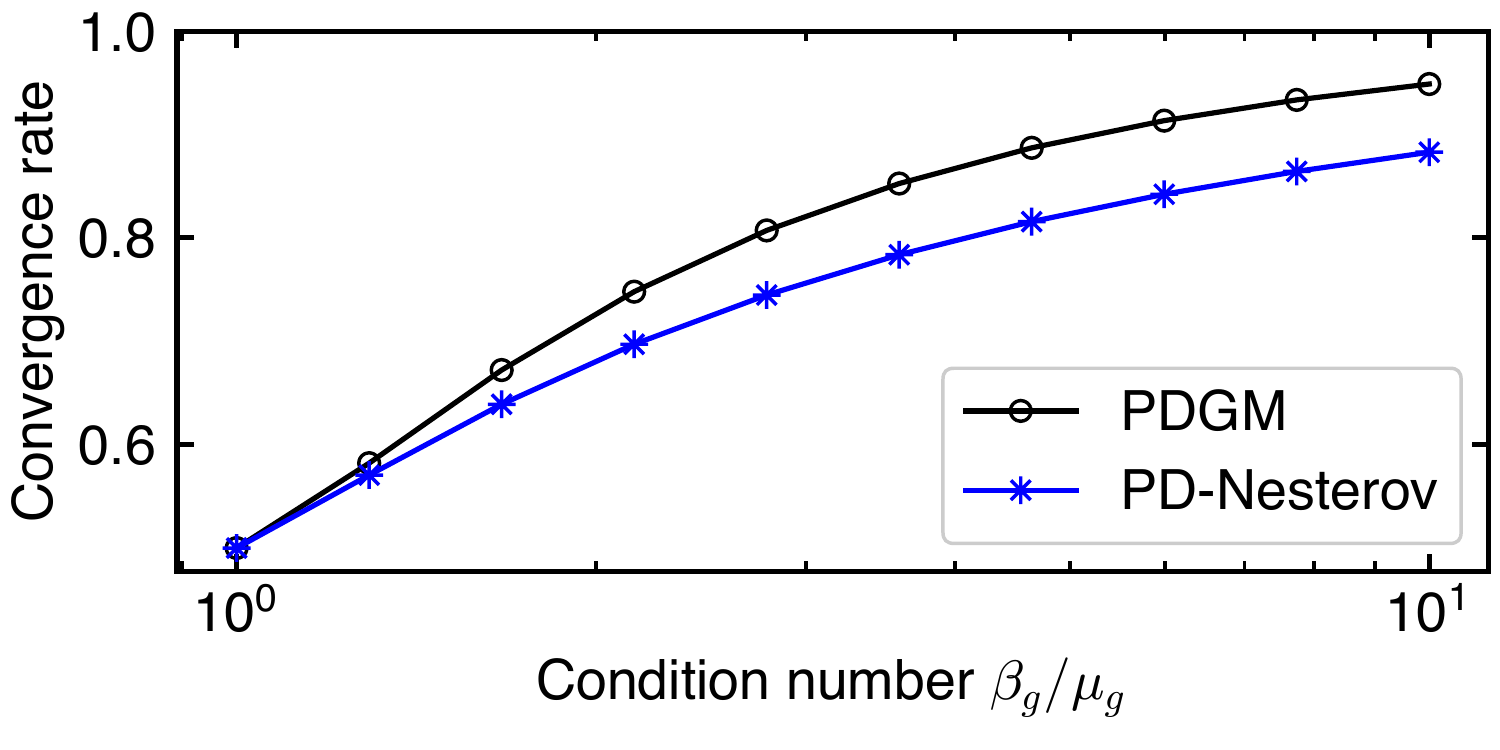}
\par\end{centering}
\caption{Convergence rate of the gradient method in~(\ref{eq:x_approx}) with different inexact oracles: gradient descent (in black, which is equivalent to the \ac{PDGM}) and Nesterov's accelerated method (in blue). \label{fig:pd-nesterov}}
\end{figure}

For a numerical comparison between the method in Corollary~\ref{cor:pd-nesterov} and the \ac{PDGM}, we considered a simple case where $f$ is linear, and $g$ is convex quadratic. For both methods, we chose $\eta_{2}=1/\betag$ and numerically searched for $\eta_{1}$ that achieved the best exponential convergence rate. Fig.~\ref{fig:pd-nesterov} shows the convergence rate for different condition numbers $\betag/\mug$. It can be seen that the method in Corollary~\ref{cor:pd-nesterov} (referred to as ``PD-Nesterov'') not only ensures convergence but also leads to a faster convergence rate compared to the \ac{PDGM}.

%% file: conclusions.tex
We have studied the convergence of inexact gradient methods in which the gradient is provided by what we refer to as a dynamic inexact oracle. When the gradient corresponds to the solution of a parametric optimization problem, dynamic inexact oracles can be realized by iterative optimization algorithms. In minimax problems, when the oracle is realized by one step of gradient descent with warm starts, the corresponding inexact gradient method recovers the \ac{PDGM}. We have shown that the interaction between the gradient method and the inexact oracle can be viewed as a feedback interconnection of two dynamical systems. Using the small-gain principle, we have derived a unified convergence analysis that only depends on a high-level description of the input-output behavior of the oracle. The convergence analysis is applicable to a range of dynamic inexact oracles that are realized by first-order methods. Furthermore, we have shown how this analysis can be used as a guideline in choosing realizations of the inexact oracle for creating new algorithms.

%% file: appendix.tex
\begin{myproof}{Lemma~\ref{lem:sec_contr}}From \cite[Thm.~2.1.15]{nesterov_introductory_2004}, we have $\cnorm{\Sstate-\eta\Sout}^{2}\leq(1-2\eta\alpha)\cnorm{\Sstate}^{2}$. The result follows from the fact $(1-2\eta\alpha)^{1/2}\leq1-\eta\alpha$. \end{myproof}

\begin{myproof}{Lemma~\ref{lem:small-gain}}Consider a single-input single-output linear system whose input $u$ and output $y$ are described by $y^{k+1}=ay^{k}+bu^{k}$, where $a\in[0,1)$ and $b\geq0$. It can be shown that the $\ell_{2}$-gain of the system is given by $b/(1-a)$. The result then follows from the (usual) small-gain theorem for feedback interconnections (see, e.g.,~\cite[Thm.~5.6]{khalil_nonlinear_2002}) and the (discrete-time) comparison lemma (see, e.g.,~\cite[Thm.~1.9.1]{lakshmikantham_stability_2015}).\end{myproof}

\begin{myproof}{Theorem~\ref{thm:io_convergence}}Define $\bar{\xi}_{i}^{k}\deq\xi_{i}^{k}-\nabla g^{*}(Ax^{k})$ ($i=1,2$), $\bar{v}^{k}\deq v^{k}-\nabla g^{*}(Ax^{k})$, and $e^{k}\deq y^{k}-\nabla g^{*}(Ax^{k})$. In the new variables, the dynamics~(\ref{eq:x_approx}) can be rewritten as~(\ref{eq:x_approx_e}), and the dynamics~(\ref{eq:dio_1st}) of $\mathcal{G}\io$ can be rewritten as 
\begin{align*}
\bar{\xi}^{k+1} & =A\io\bar{\xi}^{k}+B\io\nabla_{1}\tilde{g}(v^{k},x^{k})+\Bx\left[\phi(x^{k+1})-\phi(x^{k})\right].\\
\bar{v}^{k} & =C\io\bar{\xi}^{k},\qquad e^{k}=E\io\bar{\xi}^{k},
\end{align*}
where $\phi(x^{k})=\nabla g^{*}(Ax^{k})$ and $\Bx=-[\begin{array}{cc}
I & I\end{array}]^{T}$. 

Because $g\in\mathcal{S}(\mug,\betag)$, using Proposition~\ref{prop:smooth_and_scvx}, we have $(\bar{v}^{k},\nabla_{1}\tilde{g}(v^{k},x^{k}))=(v^{k}-\nabla g^{*}(Ax^{k}),\nabla g(v^{k})-Ax^{k})\in\sect(\mug,\betag)$. Since Assumption~\ref{assm:ABC_io} holds, we have
\begin{equation}
\cnorm{\bar{\xi}^{k+1}}_{P}\leq\rho_{2}\cnorm{\bar{\xi}^{k}}_{P}+\dxtoxi\cnorm{x^{k+1}-x^{k}}\label{eq:xi_bound_dx}
\end{equation}
for some $P\succ0$, $\rho_{2}\in[0,1)$, and $\dxtoxi>0$. The existence of $\dxtoxi$ is ensured by the Lipschitz continuity of $\phi$ and the equivalence of norms in finite dimensions. 

The second term on the right side of~(\ref{eq:xi_bound_dx}) can be further bounded by making use of~(\ref{eq:x_approx_e}). Let $\hat{x}^{k}\deq x^{k}-x\opt$, we have
\begin{align*}
\cnorm{x^{k+1}-x^{k}} & =\eta_{1}\cnorm{\nabla p(x^{k})+A^{T}e^{k}}\\
 & =\eta_{1}\cnorm{\nabla p(x^{k})+A^{T}E\io\bar{\xi}^{k}}\\
 & \leq\eta_{1}(\beta_{p}\cnorm{\hat{x}^{k}}+\gainAE\cnorm{\bar{\xi}^{k}}_{P})
\end{align*}
for some $\gainAE>0$, where we have used the equivalence of norms again. Substituting this into~(\ref{eq:xi_bound_dx}), we have 
\begin{equation}
\cnorm{\bar{\xi}^{k+1}}_{P}\leq\eta_{1}\dxtoxi\betap\cnorm{\hat{x}^{k}}+(\rho_{2}+\eta_{1}\gainAE)\cnorm{\bar{\xi}^{k}}_{P}.\label{eq:sgl_io_xi}
\end{equation}
In the meantime, because the $x$-update~(\ref{eq:x_approx_e}) is given by the gradient method, when $\eta_{1}\in(0,2/(\mup+\betap)]$, from Lemma~\ref{lem:sec_contr}, we have 
\begin{equation}
\cnorm{\hat{x}^{k+1}}\leq\rho_{1}\cnorm{\hat{x}^{k}}+\eta_{1}\gainAE\cnorm{\bar{\xi}^{k}}_{P},\label{eq:sgl_io_x}
\end{equation}
where $\rho_{1}=1-\alphap\eta_{1}$ for $\alpha_{p}$ defined in Theorem~\ref{thm:pdgm_convergence}. 

Apply the small-gain lemma (Lemma~\ref{lem:small-gain}) to~(\ref{eq:sgl_io_xi}) and~(\ref{eq:sgl_io_x}). In order to ensure convergence, we need
\begin{equation}
\begin{gathered}\rho_{1}=1-\alphap\eta_{1}<1,\qquad\rho_{2}+\eta_{1}\gainAE<1\\
\eta_{1}\gainAE\cdot\eta_{1}\dxtoxi\betap<(1-\rho_{1})(1-\rho_{2}-\eta_{1}\gainAE).
\end{gathered}
\label{eq:sgl_general}
\end{equation}
A straightforward algebraic manipulation shows that the last condition in~(\ref{eq:sgl_general}) is equivalent to $\eta_{1}<\alphap(1-\rho_{2})/(\gainAE(\alphap+\dxtoxi\betap))$. Therefore, when $\eta_{1}$ is small enough and strictly positive, all the conditions in~(\ref{eq:sgl_general}) are satisfied, which implies that the joint recursion consisting of~(\ref{eq:x_approx_e}) and~(\ref{eq:dio_1st}) converges exponentially.\end{myproof}